\documentstyle[11pt,amssymb,amscd]{article}

\newtheorem{theorem}{Theorem}[section]
\newtheorem{prop}[theorem]{Proposition}
\newtheorem{corr}[theorem]{Corollary}
\newtheorem{lemma}[theorem]{Lemma}
\newcommand{\F}{{\cal F}}
\newcommand{\E}{{\cal E}}

\newcommand{\Q}{Quot}
\newcommand{\Or}{{\cal O}^{\oplus r}}
\newcommand{\Oo}{{\cal O}}
\newcommand{\Hi}{Hilb}
\newcommand{\C}{{\Bbb C}}

\begin{document}

\centerline{\Large{\bf Moduli of Sheaves on Surfaces and}}
\centerline{\Large {\bf Action of the Oscillator Algebra}}

\bigskip

\centerline{Vladimir BARANOVSKY\footnote{email: barashek@math.uchicago.edu}}

\bigskip

\centerline{Department of Mathematics, University of Chicago}
\centerline{5734 S.University ave., Chicago, IL 60637, USA}

\bigskip

\section*{Introduction}

Let $S$ be a smooth complex  projective surface and $Hilb^n(S)$ the
Hilbert  scheme of all length $n$  zero-dimensional subschemes of $S$. 
It is known (cf. [Fo]) that $Hilb^n(S)$ is a smooth projective variety
of dimension $2n$. The structure of the cohomology ring of $Hilb^n(S)$
for a fixed  $n$ is rather  difficult to understand.  However, when we
consider  the direct sum   $\bigoplus_{n \geq 0} H^*(Hilb^n(S))$  (all
cohomology  in this  paper will  be   with complex  coefficients) the
picture becomes more comprehensible.

Firstly, for any complex smooth algebraic variety $X$ of dimension
$d$, let $P_t(X)$ be
the {\it  shifted} Poincar\'e  polynomial  $ \sum_{i=0}^{d}
t^{i-d} \cdot\dim_{\C}  H^i(X)$. It  was shown by
G\"ottsche [G\"o1] that, for any smooth quasi-projective surface $S$
$$
\sum_{n \geq 0} q^n P_t(Hilb^n(S)) = \prod_{l=1}^{\infty}
{(1+t^{-1}q^l)^{b_1(S)} (1 + tq^l)^{b_3(S)} \over 
(1- t^{-2}q^l)^{b_0(S)} (1 - q^l)^{b_2(S)} (1 -
t^2q^l)^{b_4 (S)}},
$$
where $b_i(S)$ are the Betti numbers of $S$.

Vafa and  Witten [VW] have noticed   that the right  hand  side of the
formula above   is an irreducible  character   of the  oscillator  (or
Heisenberg/Clif\-ford)  algebra   ${\cal  H}$  defined  for   a smooth
projective surface $S$ as follows:

(a) ${\cal H}$ is generated by elements $p^{\alpha}_i$, 
$\alpha \in H^*(S)$, $i \in {\Bbb Z} \setminus 0$;

(b) $[p^{\alpha}_i, p^{\beta}_j] = \langle \alpha, \beta \rangle
i \cdot \delta_{i+j, 0}$;

\noindent where the commutator above is to be understood in 
the graded sense and $\langle \cdot, \cdot \rangle$ is the 
intersection form on $S$.  The action of the oscillator algebra 
can be interpreted as follows:
there is a standard realization of Fock representation in 
 the space of all symmetric polynomials in infinitely 
many variables. This realization appears from
the cohomology of symmetric powers $Sym^n(S)$ which are related to 
Hilbert schemes of points via the Hilbert-Chow morphism $Hilb^n(S)
\to Sym^n(S)$.

\bigskip

And indeed, Nakajima and Grojnowski (cf. [Na], [Gr]) have constructed
the expected action of ${\cal H}$ on $\bigoplus_{n \geq 0} H^*(Hilb^n(S),
\C)$ using some explicit cycles in the products of  Hilbert schemes.

\bigskip

  The Hilbert scheme can be viewed as a moduli space of rank one 
sheaves on $S$ with trivial determinant and $c_2 = n$.  It was
conjectured by Vafa and Witten in [VW] (see also [Na], [Gr])
that there should be a higher rank extension of the results above.
This paper provides such a generalization  (cf. Theorem 4.1). 
 The key obsevation which
stands behind our arguments is that in certain cases a non-compact
moduli space ${\cal M}$ may admit two different natural 
compactifications ${\cal M}_1$, ${\cal M}_2$ and a map ${\cal M}_1 
\to {\cal M}_2$ which  
(partially) resolves the singularities of ${\cal M}_2$ and restricts 
to identity
on the copies of ${\cal M}$ in ${\cal M}_1$ and ${\cal M}_2$. 

The  list of examples includes moduli of maps of curves to flag 
varieties (${\cal M}_1$  and ${\cal M}_2$ are Laumon and 
Drinfeld compactifications, respectively, cf [Ku]), instantons 
on ALE spaces (quiver variety and the Uhlenbeck-type compactification
cf. [N]), moduli of abelian varieties (Voronoi and Satake
compactifications) and, finally, moduli of stable bundles on surfaces
(Gieseker and Uhlenbeck compactifications).
 Thus, the Hilbert scheme is replaced  by  the {\it Gieseker
  moduli space} $M^G(r, n)$ of stable  torsion-free sheaves and the role
of the symmetric power is played by the Uhlenbeck compactification
$M^U(r, n)$.  

  The first key result of this paper is that the fibers of the
natural map $M^G(r, n) \to M^U(r, n)$ are irreducible of expected
dimension (this generalization of the results by Brian\c con and 
Iarrobino (cf. [Br], [Ia]) was also proved independently by Ellingsrud 
and Lehn, cf. [EL]).  This leads to the second important theorem 
saying that, in the coprime case, the natural morphism $M^G(r, n) 
\to M^U(r, n)$ is strictly semismall in the sense of Goresky-MacPherson 
(if a certain technical condition
is satisfied: the surface $S$ is allowed, for
instance, to be rational, birationally ruled, K3 or abelian). This
result was originally conjectured by V. Ginzburg and the corresponding
statement is true for some of the moduli spaces mentioned above.
We proceed further to give a natural generalization of the rank one
correspondences  yielding  an action   of   the oscillator algebra  on
$\bigoplus_{n  \geq   0}  H^*(M^G(r,  n))$. We    show  how to  extend
Nakajima's  proof of commutation relations to the higher ranks case.
Some parts of the proof require a more detailed analysis of the geometry
of moduli spaces, than in the case of Hilbert schemes (see Sections 5 
and 6).

  It was communicated to us by G. Moore that $L^2$-cohomology
of Uhlenbeck compactifications have recently appeared in string theory. 
Since $L^2$-cohomology coincides in 
several interesting cases with intersection homology, this provides
another motivation for studying these homology groups.
This paper can be regarded as a first step in application of
intersection homology to moduli spaces of sheaves on surfaces.

\bigskip

{\bf Acknowledgements.} The author thanks his  advisor V. Ginzburg for
stating the problem and  for his valuable  comments and  support.  The
author is also grateful to L. G\"ottsche, A.S. Str\o mme  and Z. Qin for
their useful advices and remarks. Finally, the author thanks 
Max-Planck Institut f\"ur Mathematik where this work was carried out for
its hospitality and excellent research conditions.

\tableofcontents

\section{Moduli spaces.}

Let $S$ be a smooth complex projective surface. We denote by $Sym^n S$
the $n$-th symmetric power of $S$. Choose and fix an ample line bundle
$H$, an arbitrary line bundle $L$ and a positive integer $r$.

Consider  the moduli   space  $N(r,n)$ of Gieseker  $H$-stable  vector
bundles $E$ of  rank $r$ with  fixed determinant $L$  and $c_2(E) = n$
(cf. [Ma]). We  will only consider the  case when $gcd(r, c_1(L) \cdot
c_1(H)) = 1$. Since the bundles $L$ and $H$ will be fixed we drop them
from notation.

The moduli space $N(r, n)$ is non-compact and it can be compactified in 
two different ways. The  first   compactification is the    {\it
  Gieseker moduli space}  $M^G(r,  n)$ of all $H$-stable  torsion-free
sheaves $E$ of rank $r$ with $det(E) = L$ and $c_2(E) = n$ (cf. [Ma]).
The  second compactification  is  the   {\it Uhlenbeck  moduli  space}
$M^U(r, n)$ which can be described as follows. Take the disjoint union
$\coprod_{s=0}^{\infty} N(r,  n-s)  \times Sym^s S$.  Then Uhlenbeck's
theorem [Uh] says that any  sequence of points in   one piece of  this
disjoint union  has a subsequence that  converges (in some sense) to a
point in another piece.  This endows $\bigcup_{s=0}^{\infty} N(r, n-s)
\times Sym^s S$  with a topology of a  compact space.   Following [DK]
one can show  that this topology is  even metrizable. This topological
space is called the Uhlenbeck  compactification  $M^U(r, n)$ of  $N(r,
n)$.

Note that the union above is in fact finite since a necessary condition
for non-emptiness of $N(r, k)$ is the Bogomolov inequality $2 r k -
(r-1) (c_1(L))^2 > 0$. We can always tensor all our bundles and sheaves
with a high power of $H$ and achieve $c_1^2(L) > 0$. Then Bogomolov
inequality implies that $k=n-s$ should be at least positive.

Hence without loss of generality we can assume that $0 \leq s < n$.

\medskip

In general the  two compactifications above may  be quite difficult to
investigate: $N(r, n)$ may not be dense in $M^G(r, n)$ or $M^U(r,
n)$, some components may be of dimension higher than expected,
etc.  For this reason we introduce a

\bigskip

{\bf Technical Condition.} 

\medskip

\noindent
(a) The integers $r$ and $d:= c_1(L) \cdot c_1(H)$ are coprime.

\noindent
(b) Either the canonical bundle $K_S$ is trivial or $c_1(K_S) \cdot
c_1(H) < 0$.

\bigskip 

The main reason for imposing (a) and (b) is that
they ensure (cf. [HL]) that for any $n \geq  0$ the moduli space 
$M^G(r, n)$ is either empty or smooth of expected dimension 
\begin{equation}
\label{condition}
\dim M^G(r, n) = 2rn  - (r-1) (c_1(L))^2 - (r^2 - 1) \chi({\cal   O}_S)  
+ h^1({\cal O}_S).
\end{equation}

\smallskip

This condition  (b)  is automatic  when   $(-K)$ is
represented by an effective curve. Note that since $N(r, n)$ is an
open subset of $M^G(r, n)$ ({\it loc. cit.}), it also has to be smooth. 

Since  $r$ and  $d$  are assumed to  be coprime, the two possible
notions of stability (Mumford-Takemoto and Gieseker) coincide.
\medskip

{\bf Examples.} If we require that  $K_S$ is trivial or effective then
$S$ can  be a Del  Pezzo surface, K3 or abelian  surface. If we make a
special choice of $H$ as above the list of examples will extend to all
rational surfaces and birationally ruled surfaces.

It is  known (e.g.  [DK 10.3.4])  that for some  elliptic  surfaces or
surfaces of general type the condition (\ref{condition}) may fail.

\section{Gieseker to Uhlenbeck.}

The Gieseker and  Uhlenbeck compacifications are higher rank analogues
of the Hilbert scheme of all length $n$ subschemes $Hilb^n(S)$ and the
symmetric power  $Sym^n(S)$, respectively. Here  we prefer to think of
$Hilb^n(S)$ as    parametrizing  rather   the ideal  sheaves,    i.e.  
torsion-free sheaves of  rank one with  $c_2=n$ and trivial
determinant. One has a natural map
$a_1: Hilb^n(S) \to Sym^n(S)$:
$$
a_1: J_{\xi} \mapsto \sum_{x_i \in Supp\; \xi} mult_{x_i}(\xi) \cdot  x_i
$$
where $J_{\xi}$ is an ideal sheaf of a subscheme $\xi$ and the formal
sum above is viewed as an element of $Sym^n(S)$.

There exists a natural map $a_r: M^G(r, n) \to M^U(r, n)$ generalizing
$a_1$. To define $a_r$ first note that, for any torsion-free sheaf $\F$,
the double-dual $\F^{**}$ is reflexive. It is known (cf. [OSS]) that any
reflexive sheaf on a smooth variety is locally free on the complement of
a closed subset of codimension $\geq 3$. Hence in our case ($\dim S=2$)
the double-dual is necessarily locally free.

Consider the short exact sequence: $$ 0 \to \F \to \F^{**} \to A_{\F} \to
0. $$ The quotient sheaf $A_{\F} = \F^{**}/\F$ is supported on finitely
many points (i.e. $A_{\F}$ is an Artin sheaf). Denote by $l_{\F}$ the
length of $A_{\F}$. We define $a_r(\F)$ by $$ a_r: \F \mapsto (\F^{**},
\sum_{x_i \in Supp\; A_{\F}} mult_{x_i}(A_F) \cdot x_i) $$ where the
image is a point in $ N(r, n-l_{\F}) \times Sym^{l_{\F}}(S) \subset
M^U(r, n)$.

\medskip

{\bf Remark.} Another consequence of  $r$ and $d= c_1(H) \cdot c_1(L)$
being  {\it coprime}  is   that any  semistable  sheaf  is necessarily
stable.  This  ensures that $\F^{**}$  is stable  if  $\F$ is, hence the
definition above makes sense.   If $r$ and  $d$  are not coprime,  the
moduli space $M^G(r, n)$ parametrizes  only S-equivalence classes (cf. 
[Ma]) of  semistable sheaves and  one has to be  more  careful to make
sure that $a_r$ is well-defined on such classes (cf. [Li]).

\medskip 

It takes quite a  lot of technical work  (see [Li], [Mo]) to show that
the   map   $a_r$ is  continuous (this  is    because one has   to use
gauge-theoretic Uhlenbeck Compactness Theorem and consider $M^U(r, n)$
as a topological  space). We will  not attempt to repeat  the argument
here  referring the   interested   reader to  the    references cited. 
Moreover, from now on we {\it assume that  $M^U(r, n)$ has a structure
of a projective algebraic variety  and $a_r$ is algebraic}.  In rank
2 case this  was proved in [Li] and  we intend to give  an alternative
construction for the general case in a forthcoming paper.

We need a  description of the fiber  of $a_r$ over  an arbitrary point
$p=(E, \sum_i   m_i x_i)$ in   $M^U(r, n)$.  Here  $m_i$  are positive
integers   and the  points $x_i$  are  pairwise   distinct. The  fiber
$a_r^{-1}(p)$ parametrizes all the quotients $E \to A \to 0$ where the
sheaf    $A$  is supported at       the points $x_i$  with  prescribed
multiplicities $m_i$.  Since the question is local, the fiber does not
depend on $E$ and points $x_i$  but only on  $r$ and $m_i$. To be more
precise, let $Quot(r, n)$ denote the {\it punctual} Quot scheme of all
quotients ${\cal O}^{\oplus  r} \to A$ of  fixed length $n$  which are
supported at a fixed point $x$. This scheme depends only on completion
$\widehat{{\cal O}}_x$ of the  local ring ${\cal  O}_x$ of  $x$. Hence
different points $x$ lead  to isomorphic Quot  schemes and we drop $x$
from notation. The following statement is immediate

\begin{prop}
The fiber $a_r^{-1}(p)$ over the point $p = (E, \sum_i m_i x_i)$ is 
isomorphic to the product of punctual Quot  schemes $\prod_i Quot(r, m_i)$.
\hfill $\Box$
\end{prop}

 The proposition above motivates the following theorem (cf. [Ba] and
[EL]):

\begin{theorem} The punctual Quot scheme $Quot(r, n)$ is irreducible of
dimension $rn - 1$.
\end{theorem}
{\it Proof.} See Appendix. $\Box$

\bigskip

 The first application of Theorem 2.2 is to show  that, under the
technical assumption above,  the moduli spaces $M^G(r, n)$
and $M^U(r, n)$ are well-behaved.

\begin{theorem} Assume that the technical condition of Section 1 is
satisfied. Then the following statements hold

(a) $N(r, n)$ is dense in $M^G(r, n)$;

(b) $N(r, n)$ is dense in $M^U(r, n)$;

(c) Any irreducible component of $N(r, n-s) \times Sym^s(S)$ intersects
the closure of a unique component of $N(r, n)$ (and hence by (b) is
contained in it).

\end{theorem}
{\it Proof.} Note that (a) implies (b) since $a_r$ is surjective and 
one-to-one on the copies of $N(r, n)$ in $M^G(r, n)$ and $M^U(r, n)$,
respectively. 

To prove (a) suppose that $N(r, n)$ is not  dense in $M^G(r, n)$. Then
there exists a component of $M^G(r, n)$ such that  generic point of it
corresponds to a non-locally free sheaf.  This means that the image of
this component in $M^U(r, n)$ is a subset of $\bigcup_{i \geq 1} N (r,
n-s) \times   Sym^s(S)$. Since all  the components  of $N(r,  n)$ have
expected dimension (and not more than that), Theorem 2.2 above implies
that the  dimension of $a_r^{-1}(\bigcup_{i \geq 1}  N (r, n-s) \times
Sym^s(S))$ is strictly less than the expected dimension of $M^G(r, n)$
which is impossible.

Finally, if  (c) were false, some  irreducible  component $X$ of $N(r,
n-s) \times  Sym^s(S)$  would intersect the  closures  of at least two
different components of $N(r, n)$. Since $M^G(r,  n)$ is smooth, these
would  mean that   $X$  would belong  to  the  image of two  different
connected components of $M^G(r, n)$. But  this is impossible since all
fibers of $a_r$ are irreducible. \hfill $\Box$

\begin{corr}  $M^U(r, n)$ is a disjoint union of the closures of
irreducible components of $N(r, n)$. $\Box$
\end{corr}

\section{Stratifications and semi-small maps.} 

  From now on we will assume that the technical condition on the 
pair $(S, H)$ is satisfied.   

Recall briefly the  results on symmetric  products  and  Hilbert
schemes.
  
Let $\mu$  be a partition  of $n$. Any such   $\mu$ can be represented
either by a non-increasing sequence $\mu=(\mu_1 \geq \mu_2 \geq \ldots
\geq \mu_m > 0)$ with $\sum \mu_i =  n$ or in the form $1^{m_1}2^{m_2}
\ldots n^{m_n}$, where  $m_i$ is the number  of parts $\mu$ which  are
equal to $i$ (hence  $m_i \geq 0$).  Note that  $m_1 + m_2 +  \ldots +
m_n$ is equal to the number $m$ of non-zero parts of $\mu$.

 The symmetric product $Sym^n(S)$ has a natural stratification by locally
closed strata  labeled by partitions of $n$.

The stratum $Sym^n_{\mu}(S)$ is the set of formal sums of the type 
$\sum \mu_i x_i$ with $x_i$ pairwise distinct (veiwed as elements of 
$Sym^n S$). Note that $Sym^n_{(1, \ldots, 1)}$ is
a dense open subset of $Sym^n(S)$ and in general $Sym^n_{\mu}(S)$ is
isomorphic to a dense open subset of $Sym^{m_1}(S) \times \ldots \times
Sym^{m_n}(S)$. A generic element $(y_1, \ldots, y_n)$ in the latter
product corresponds to the point $y_1 + 2y_2 + \ldots + ny_n$ in
$Sym^n_{\mu}(S)$.

  Let $\pi: Z \to Y$ be a proper projective morphism  of algebraic 
varieties. Suppose that $Y$ decomposes into a finite number of locally
closed strata: $Y = \bigcup_{\mu} Y_{\mu}$ and choose an arbitrary
point $y_{\mu} \in Y_{\mu}$. Assume that the restriction 
$\pi: \pi^{-1}(Y_{\mu}) \to Y_{\mu}$ is a topological fiber bundle
with fiber $\pi^{-1}(y_{\mu})$.

\bigskip

{\bf Definition.}  (cf.  [BM] or [CG,   Chapter 8]) The map  $\pi$  is
called strictly semi-small if it satisfies
$$
2\; dim \; \pi^{-1}(y_{\mu}) = codim\; Y_{\mu}
$$
for any stratum $Y_{\mu}$.

The  following proposition is an immediate   consequence of results of
Brian\c{c}on and Iarrobino (cf. [Br], [Ia]).

\begin{prop} (cf. [GS])

  The morphism $a_1: Hilb^n(S) \to Sym^n(S)$ is strictly semi-small,
with respect to the stratification given by $Sym^n_{\mu}(S)$.
\end{prop}

We  will  give a generalization of  this  proposition  to the  case of
arbitrary rank. Exactly as in [GS] this will lead to some formulas for
Poincar\'e polynomials.  These  formulas  will be later   interpreted in
representation-theoretic terms.

The  first step  is to  stratify $M^U(r,  n)$ appropriately since  the
natural strata coming from the definition of $M^U(r,  n)$ are too big. 
The new finer strata will be labeled by pairs $(s, \mu)$  where $s < n
$ is  a non-negative integer such  that $N(r, n-s)$  is non-empty, and
$\mu =  (\mu_1  \geq  \mu_2 \geq  \ldots  \geq  \mu_m \geq   0)$  is a
partition   of  $s$. For  each  pair $(s,   \mu)$ we consider $M^U_{s,
  \mu}(r, n)  = N(r, n-s) \times  Sym^n_{\mu}(S)$ which is naturally a
locally closed subset of $M^U(r, n)$.  Of course, this is nothing but
  a "common refinement" of the  natural partition of $M^U(r, n)$ and
the above partition above of  the symmetric product. 

Part (b)  of
the   next  theorem was  originally   conjectured by  V.  Ginzburg. It
provides   a starting   point for   our  generalization of  Nakajima's
construction.

\begin{prop} Let $s$ and $\mu = (\mu_1 \geq \mu_2 \geq \ldots 
  \geq \mu_m >  0)$ be as  above and let  $M^U_{s, \mu}(r,  n)$ be the
  stratum associated to the pair $(s, \mu)$. Then
  
  (a) For any point $x_{s, \mu} \in M^U_{s,  \mu}(r, n)$ the dimension
  of the fiber $a_r^{-1}(x_{s, \mu})$ is equal to $(rs - m)$ where $m$
  is the number of non-zero parts of $\mu$.

(b) The morphism $M^G(r, n) \to M^U(r, n)$ is strictly semi-small with
respect ot the stratification given by $M^U_{s, \mu}(r, n)$.
\end{prop}
{\it Proof.}   To prove (a) note that by Proposition 1 
the fiber is isomorphic to 
$\prod_{i=1}^m Quot(r, m_i)$ and by Theorem 2
its dimension is equal to $\sum_{i = 1}^{m} (r\mu_i -1) = r(\sum_{i =1}^m 
\mu_i) - m = rs - m$.

To  prove (b) we have  to  show that  $codim\; M^U_{s,  \mu} = 2 dim\;
a_r^{-1}(x_{s, \mu}) = 2(rs -m)$. In fact, $codim\; M^U_{s, \mu}(r, n)
= dim\;N(r, n) - dim\; N (r, n-s) - dim\;  Sym^n_{\mu}(S) = 2rs - 2m =
2(rs-m)$. $\Box$

\bigskip 

  The semi-smallness result will allow us to relate homological 
invariants of $M^G(r, n)$ and $M^U(r, n)$. Recall (cf. [BBD]) that for 
any algebraic variety $X$, there exists a remarkable complex of sheaves  
$IC_X$ ({\it intersection cohomology}  complex) such that its 
cohomology groups $IH^*(X) = H^*(X, IC_X)$ ({\it intersection} or 
{\it Goresky-MacPherson} homology) satisfy Poincar\'e duality.
If $X$ is a smooth algebraic variety (such as $M^G(r, n)$) or a quotient
of a smooth variety by a finite group action (such as $Sym^n(S)$) then
up to shift $IC_X$ is just the constant sheaf ${\Bbb C}$ (cf. [BBD]
of [GS]).

 Following [Sa, 1.13] one can define pure Hodge structure on the
intersection homology of any complex algebraic variety $X$.

 We will need a simplified version of Borho-MacPherson formula 
for the direct image of the intersection homology complex under a
projective semismall morphism. This formula is a direct application of
 the Decomposition Theorem due to Beilinson-Berstein-Deligne-Gabber
(cf. [BBD]).

\begin{prop}  Suppose that a projective morphism $Z \to Y$ of algebraic
varieties is strictly semi-small with respect to some stratification 
$Y = \bigcup Y_{\mu}$. Suppose further that for any point $y \in Y$ the
fiber $\pi^{-1}(y)$ is irreducible. Then
$$
\pi_*(IC_Z) = \bigoplus_{\mu} IC_{\overline{Y}_{\mu}}. \quad \quad \Box
$$ 
\end{prop}

\begin{corr} One has the following direct sum decomposition in 
the derived category of complexes of sheaves:
$$
(a_r)_* IC_{M^G(r, n)} = \bigoplus_{s \in \{0, \ldots n\}, \,\mu \in P(s) } 
IC_{\overline{M}^U_{s , \mu }(r, n)}
$$
where $P(s)$ is the set of all partitions of $s$ and 
$IC (\overline{M}^U_{s , \mu })$ denotes the $IC$-complex
on the closure of the stratum $M^U_{s, \mu} \subset M^U(r, n)$. $\Box$
\end{corr}

We intend to use the  formula above by  taking the cohomology of  both
sides.  The $IC$-complexes supported on  the  closures of the  smaller
strata can be understood with the help of the following

\begin{prop}
  Let  $(s, \mu)$ and $m_i$ be as above and denote by $Sym^{\mu}(S)$
the direct product $ Sym^{m_1}(S) \times \ldots \times Sym^{m_n}(S)$. 
Then there exists a finite birational
morphism respecting the induced stratifications:
$$
\pi_{n, s, \mu}: M^U(r, n-s) \times Sym^{\mu}(S) \to \
\overline{M}^U_{s , \mu }(r, n)
$$
\end{prop}
{\it Proof.} Left as an exercise to the reader. $\Box$

\medskip 

Now [GS, Lemma 1] implies that in this situation one can deduce a

\begin{corr} In the notation of the previous proposition, one has
the following equality in the derived category of sheaves:
$$
(\pi_{n, s, \mu})_* (IC_{M^U(r, n-s) \times Sym^{\mu}(S)}) = 
IC_{\overline{M}^U_{s, \mu}(r, n)}. \quad \quad \Box
$$
\end{corr}
\medskip

Recall that intersection homology complex of a space $X$ is defined in
such  a way that  it has non-trivial  (global)  hypercohomology in the
range between $(-k)$ and $k$ where $2k$ is the real dimension of $X$.

\medskip

{\bf Definition.} The shifted intersection   homology  Poincar\'e
polynomial $P_t(X)$ is defined by the formula
$$
P_t(X):=  \sum_{-\dim_{\C}   X}^{\dim_{\C} X} (\dim_{\C}  H^{i}(X,
IC_X))\cdot t^i
$$
Note that for a   smooth $X$ this is   just  the usual Poincar\'e
polynomial multiplied by $t^{-\dim_{\C} X}$.

Similarly, one defines a shifted Hodge polynomial $P_{x,y}(X)$ using the
pure Hodge structure on $IH^*(X)$. When $X$ is smooth this coincides
with the usual Hodge polynomial multiplied by $(xy)^{- \dim_{\C} X/ 2}$.
However, in general it will {\it not} be a shift of the virtual Hodge
polynomial of $X$.

\bigskip

With this preparation, we can deduce the first main result of the paper.

\begin{theorem} 
One has the following identity between Poincar\'e polynomials of
$M^G$ and $M^U$:
$$
{\sum_{n   =0}^{\infty} q^n   P_t   (M^G(r, n))  \over \sum_{n  =
  0}^{\infty}    q^n  P_t  (M^U(r,      n)) }  =     
\prod_{l=1}^{\infty} {    (1   +  t^{-1}q^l)^{b_1(S)}  (1     +  t^{1}
  q^l)^{b_3(S)} \over (1 - t^{-2} q^l)^{b_0(S)}  (1 - q^l)^{b_2(S)} (1
  -t^{2} q^l)^{b_4 (S)} } 
$$
In short,  the series for  $M^G$ is obtained from  the series
for $M^U$ when multiplied by the G\"ottsche's formula.
\end{theorem}
{\it Proof.} Following Nakajima's notations,  denote by $a_m(t)$ the 
(shifted) Poincar\'e polynomial $P_t(Sym^m(S))$.
Then a formula due to MacDonald [Mc1] says that 
$$
\sum_{m=0}^{\infty} q^m a_m(t) = { (1 + t^{-1}q)^{b_1(S)}(1 + tq)^{b_3(S)} \over
(1-t^{-2}q)^{b_0(X)}(1-q)^{b_2(S)} (1-t^2q)^{b_4(S)}}.
$$

 Taking the cohomology of both sides in the formula of Corollary 3.4 and
using Corollary 3.6 one obtains:
$$
\sum_{k =0}^{\infty} q^k P_t (M^G(r, k)) = \sum_{k=0}^{\infty}
\sum_{s=0}^k \sum_{\mu \in P(s)} q^k P_t (M^U(r, k-s)) \cdot
P_t(Sym^{\mu} S)=
$$

$$
= \sum_{k =0}^{\infty} \sum_{{ 0 \leq s \leq k \atop \mu \in
  P(s)}} q^{k-s} P_t (M^U(r, k-s))
\cdot  q^s  P_t (Sym^{\mu} S)=
$$

$$
= 
\Big(\sum_{n = 0}^{\infty} q^n P_t (M^U(r, n)) \Big) \cdot
\Big(\sum_{{s \geq 0 \atop \mu \in P(s)}} q^s P_t (Sym^{\mu} S) \Big).
$$

Represent all partitions $\mu$ in the formula above
as $(1^{m_1}, 2^{m_2}, \ldots, s^{m_s})$.  Then

$$
\sum_{{s \geq 0 \atop \mu \in P(s)}} q^s P_t (Sym^{\mu} S) =
\sum_{{s \geq 0 \atop \mu \in P(s)}} a_{m_1}(t) (q)^{m_1} 
 a_{m_2}(t) (q^2)^{m_2}  \ldots 
 a_{m_s}(t) (q^s)^{m_s} = 
$$

$$
= \prod_{l=1}^{\infty} \Big(\sum_{m=0}^{\infty} a_m(t)(q^l)^m \Big)
= \prod_{l=1}^{\infty} { (1 + t^{-1}q^l)^{b_1(S)} (1 + t^{1}
  q^l)^{b_3(S)} \over (1 - t^{-2} q^l)^{b_0(S)} (1 - 
  q^l)^{b_2(S)} (1 -t^{2} q^l)^{b_4 (S)} }. \; \Box
$$

\bigskip

{\bf Remarks.} 

(1) A very similar formula was discovered by G\"ottsche in [G\"o22]. 
He computes the ratio of generating functions for $M^G(r, n)$ and
$N(r, n)$ and obtains a similar product (involving more factors).

(2)  In [LQ]  Li   and Qin derive   a  formula relating  cohomology of
Gieseker spaces for   $S$  and its blowup  $\widetilde{S}$.   Applying
Theorem 3.7 we get  a  relation between  the intersection homology  of
Uhlenbeck compactifications for  $S$ and $\widetilde{S}$. This formula
is different form the relation  between usual homology of $M^U(r,  n)$
for $S$, $\widetilde{S}$ (also  found by Li  and Qin), since Uhlenbeck
compactifications are quite singular hence  their intersection homology is
different from the ordinary (co)homology.

\bigskip

{\bf Example.} When $S = {\Bbb P}^2$ and $r=2$ the generating function
for $M^G(r, n)$ was computed by Yoshioka in [Y].
In this case 
$$
\sum_{i=0}^{\infty} q^n P_t(M^G(r, n)) = \hspace{5cm}
$$
$$
\hspace{2cm}=
{\sum_{b \in {\Bbb Z}} {t^{-2b}q^{b^2} \over 1 - t^4 q^{2b-1}} 
\over \sum_{n \in {\Bbb Z}} t^{- 2n}q^{n^2}} \cdot {1 \over t^4(t^2 -1)
\prod_{l= 1}^{\infty} (1 - t^{-2}q^l)^2(1 - q^l)^2 
(1 -  t^{2}q^l)^2}
$$
Let $q = e^{2\pi i \tau}$ and $t = e^{2 \pi i z}$. 
Recall the classical theta functions:
$$
\theta_{\mu, \nu}(\tau, z) := \sum_{n \in {\Bbb Z}} (-1)^{n\nu}
q^{(n+\mu/2)^2/2} t^{n + \mu/2} \quad (\mu, \nu \in \{0, 1\}).
$$
Then by the product formula
$$
\theta_{1, 1}(\tau, z) = q^{1/8} (t^{1/2} - t^{-1/2}) \prod_{l> 0}
(1-t^{-1}q^l)(1-q^l) (1 -tq^l)
$$
for theta functions one can write
$$
\sum_{i=0}^{\infty} q^n P_t(M^G(r, n)) = 
\sum_{b \in {\Bbb Z}} {t^{-2b}q^{b^2} \over 1 - t^4 q^{2b-1}} \cdot 
{ q^{1/4} (t - t^{-1})  \over t^5\; \theta_{0, 0}(2\tau, 2z) \theta_{1,
    1}(\tau, 2z)^2}
$$

By the Theorem above we immediately deduce that 
$$
\sum_{i=0}^{\infty} q^n P_t(M^U(r, n)) = 
\sum_{b \in {\Bbb Z}} {t^{-2b}q^{b^2} \over 1 - t^4 q^{2b-1}} \cdot 
{ q^{1/8} \over t^5 \; \theta_{0, 0}(2\tau, 2z) \theta_{1, 1}(\tau, 2z)}
$$ 

{\bf  Remark.} Conjecturally,  the theta functions in the generating
series above appear  from extra  symmetries  of cohomology of the
moduli spaces.

\bigskip
 Using Saito's theory of mixed Hodge modules (cf. [Sa]) and the formula
for Hodge polynomials of symmetric powers (cf. [Bu]) one can apply the 
arguments above to prove the following 

\begin{prop}

(a) There exists an isomorpshism of pure Hodge structures:
$$
IH^*(M^G(r, n)) \simeq \bigoplus_{{s \geq 0 \atop \mu \in P(s)}} 
IH^*(M^U(r, n-s))
\otimes IH^*(Sym^{\mu}(S));
$$

(b) One has the following identity between the (shifted) intersection
homology Hodge polynomials of $M^G$ and $M^U$:
$$
{\sum_{n \geq 0} P_{x, y}(M^G(r, n)) q^n \over
\sum_{n \geq 0} P_{x, y}(M^U(r, n)) q^n} =
\prod_{l \geq 0}\prod_{i, j = -1}^1 
(1 + (-1)^{i+j+1} x^i y^j q^l)^{(-1)^{i+j+1}h^{i+1, j+1}(S)};
$$
where $h^{r, s}(S)$ are the (usual) Hodge numbers of $S$. \quad \quad $\Box$
\end{prop}

\bigskip

{\bf Remark.} Since $S$ is compact, one of course has $b_0 = b_4 = 1$ and
$b_1 = b_3$. However, G\"ottsche-Soergel's approach is valid for Hilbert 
schemes of any quasi-projective surface. It is interesting to know if
the computaion above can be extended  to the quasi-projective case.

\section {Correspondences.}

In  this   section  we   give   a   very  natural   generalization  of
correspondences used  by Nakajima and  Grojnowski and prove that their
action on the cohomology satisfies the  correct commutation relations. 
Of  course, as one  might guess from the  formula  of Theorem 3.7, the
corresponding module over   the oscillator algebra is not  irreducible
any more. In   fact,  the space of  vacuum  vector  can  be  naturally
identified with intersection homology of Uhlenbeck compactifications.

\bigskip 

Recall [CG] that, for any  two smooth projective varieties $M_1, M_2$,
a closed  subvariety $Z \subset M_1 \times  M_2$  (more generally, any
cohomology class in  $H^*(M_1  \times  M_2$)),  defines a  map   $[Z]:
H^*(M_1) \to H^*(M_2)$   given   by $[Z](\alpha)  = (p_2)_*([Z]   \cap
(p_1)^*)$, where $p_1$  and $p_2$ are the  projections onto $M_1$  and
$M_2$, respectively. If $M_3$ is a third smooth projective variety and
$Z'  \subset  M_2 \times M_3$ then  the  composition $[Z']  \circ [Z]:
H^*(M_1) \to H^*(M_3)$  can be described  as follows. Let $p_{ij}: M_1
\times  M_2 \times M_3 \to  M_i \times M_j$  be the  projection to the
product of the  $i$-th  and $j$-th factors.  Then $[Z']  \circ [Z]$ is
given   by    the   cohomology   class   $(p_{13})_*(p_{12}^*[Z]  \cup
p_{23}^*[Z']) \in H^*(M_1 \times M_3)$.

\bigskip

We want     to   apply this   construction  by    introducing  certain
correspondences in the products of moduli spaces.

\bigskip

{\bf Definition.} For $i > 0$, define a subvariety $P_{-i} \subset 
\coprod_{n=0}^{\infty} M^G(r, n) \times M(n+i)$ as a set of all pairs 
$(\F_1, \F_2)$ such that

1) $\F_1^{**} \simeq \F_2^{**}$, 

2) $\F_2 \subset \F_1$ as subsheaves of the common double  dual, 

3) $Supp(\F_1/\F_2) = \{x\}$ for some $x \in X $.

\bigskip

Note that the second condition makes sense since the double dual is
stable  and therefore simple.

 Similarly, we define $P_i$ for $i > 0$, exchanging $\F_1$ and $\F_2$.
Let us define a morphism $\Pi: P_i \to S$ by
$$
\Pi(\F_1, \F_2) = x \quad \mbox{for } (\F_1, \F_2) \in P_i,
$$
where $x$ is the unique element of $Supp(\F_1/\F_2)$. One can give a
rigorous definition of $\Pi$ using universal sheaves on moduli spaces
and [HL, Example 4.3.6] but we will not do it here.

Denote the dimension of $M^G(r, n)$ by $2rn + a$ where $a$ is a number
depending on $S, r, L$ and $H$ but  not on $n$, cf. (1.1).  Then the 
dimension of
$P_{-i} \cap M^G(n-i) \times  M^G(n)$  is equal to  $2rn  + a -ri +  1$
(this follows from Nakajima's argument [Na 8.3] and Theorem 2.2).

\bigskip

Let   $P^{\alpha}_{\pm i} \in   \prod_n  H^*(M^G(r, n  \mp i))  \otimes
H^*(M^G(r, n))$, $\alpha \in H^*(S)$, be the  images of the cohomology
classes  $\Pi^*\alpha$    on     $P_{\pm i}$.     The   cycles
$P^{\alpha}_{\pm i}$ induce  the maps  $\bigoplus_n  H^*(M^G(r,  n)) \to
\bigoplus_n H^*(M^G(r \mp i))$  via the convolution construction given
above.  Computing dimensions we can  see  that $P^{\alpha}_{-i}$
raises cohomological degree by $2(ri-1) + deg\;\alpha$.

\bigskip

The   next theorem asserts   that $\bigoplus_{n=0}^{\infty} H^*(M^G(r,
n))$ becomes a representation space of the Heisenberg/Clifford algebra
generated by $P^{\alpha}_{i}$ and $P^{\beta}_{-i}/(-1)^{ri-1}r$. It is a
direct generalization of [Na, Theorem 8.13].

\begin{theorem} The following relations hold,
$$
[P^{\alpha}_{i}, P^{\beta}_{j}] = (-1)^{ri-1} r i \delta_{i+j, 0} \langle 
\alpha, \beta \rangle \; Id \quad \quad \quad \quad if (-1)^{deg \alpha 
\cdot \deg \beta} = 1,
$$
$$
\{P^{\alpha}_{i}, P^{\beta}_{j}\} = (-1)^{ri-1} ri \delta_{i+j, 0} \langle 
\alpha, \beta \rangle \; Id \quad otherwise,
$$
where $\langle \alpha, \beta \rangle$ is the intersection form on $S$.
\end{theorem}

  Exactly as in [Na Corollary 8.16] we deduce a

\begin{corr} The cohomology groups $\oplus_{n=0}^{\infty} H^*(M^G(r, n))$
form a tensor product of an irreducible highest weight representation
of the oscillator algebra, with a trivial representation on the
intersection homology groups of Uhlenbeck compactifications.
\end{corr}

\bigskip

{\bf Example.} Let $r=2$, $i=1$, $S= {\Bbb P}^2$ and $L \simeq H \simeq 
\Oo(1)$. Choose $\alpha$ to be the fundamental class and $\beta$ the 
class of a point $x \in {\Bbb P}^2$. Then the formula above predicts 
that $[P^{\alpha}_{1}, P^{\beta}_{-1}] = (-2) Id$. This can be seen as
follows.

 It is known that $M^G(2, 1)$ is a single point 
corresponding to the twist $T_{{\Bbb P}^2} (-1)$ of the tangent bundle 
(cf. [OSS, 3.2]). The Uhlenbeck compactification $M^U(2, 2)$ can be 
described as follows. Let $V$ be the 6-dimensional space of symmetric 
$3 \times 3$ matrices.
Then $M^U(2, 2)$ can be identified with a hypersurface $H$ in $P(V)$
correspoding to all rank $\leq 2$ symmetric matrices (cf. [OSS, 4.3]). 
Its singular locus $P \simeq {\Bbb P}^2$ corresponds to symmetric
matrices of rank one. It can be deduced from (cf. [OSS, 4.3]) that 
$M^G(2, 2)$ coincides with the blowup of $P$. 

Let $D \subset M^G(2, 2)$
be the exceptional divisor. Each one-dimensional fiber of $D \to P$
corresponds to the set of sheaves $\F$ which fit a short exact sequence
$$ 0 \to \F \to T_{{\Bbb P}^2} (-1)
\to \C_x \to 0,$$ where $x \in {\Bbb P}^2 \simeq P$. Let $l$ be such a
fiber and $[pt]$ be the generator of $H^0(M^G(2, 1)) = H^0(pt)$. Since
$P^{\alpha}_{1}([pt])=0$ and $P^{\beta}_{-1}([pt])= [l]$ we expect that
$P^{\alpha}_{1}([l]) = -2[pt]$. It follows from definitions that this is
equivalent to $[D] \cdot [l]= -2$. The last equality is true since
$H$ has an $A_2$-singularity at any point of $P$.

\bigskip 

{\it Beginning of the proof of 4.1:} Most of Nakajima's proof
can be repeated word-by-word with almost no changes. In particular, we
have the following result:

\begin{theorem} For any $v \in H^*(M^G(r, n))$, the following 
relations hold:
$$
[P^{\alpha}_i, P^{\beta}_j](v) =  c_{i, n} \delta_{i+j, 0} \langle 
\alpha, \beta \rangle \; v \quad \quad \quad \quad if (-1)^{\deg \alpha 
\cdot \deg \beta} = 1,
$$
$$
\{P^{\alpha}_i, P^{\beta}_j\} (v) = c_{i, n} \delta_{i+j, 0} \langle 
\alpha, \beta \rangle \; v \quad otherwise,
$$
where $\langle \alpha,  \beta \rangle$ is  the intersection form on
$H^*(S)$ and $c_{i, n}$  is some constant  depending {\em  apriori} on
$i$ and $n$ but not on the classes $\alpha$ and $\beta$. \quad $\Box$
\end{theorem}

The proof  is a mere repetition of  the argument presented in [Na 8.4]
(where all the dimension statements are proved using Theorem 2.2).

\begin{lemma} For each $i$, the constants $c_{i, n}$ are independent
of $n$.
\end{lemma}
{\it Proof.} Let $k \neq \pm i$. By Theorem 4.3 we have $[P^{\gamma}_{-k},
[P^{\alpha}_i, P^{\beta}_{-i}]] = 0$. Take $v \in H^*(M^G(r, n))$.
Applying the double commutator to $v$, one gets 
$$
\langle \alpha, \beta \rangle (c_{i, n} - c_{i, n+k}) P^{\gamma}_{-k} v = 0.
$$
One can choose $v, \alpha, \beta, \gamma$ in such a way that
$\langle \alpha, \beta \rangle \neq 0$ and $P_{\gamma}[-k](v) \neq 0$.
Therefore, one  has $c_{i, n} =
c_{i, n+k}$ if $k \neq \pm i$. If $i \neq \pm 1$ we can take $k=1$ and
obtain $c_{i, n} = c_{i, n+1}$. If $i = \pm 1$ we take $k  = 2, 3$ and
get $c_{1,   n} = c_{1,  n+2}  = c_{1, n+3}$   which also implies that
$c_{1, n}$ is independent of $n$. $\Box$

\bigskip

In view of this lemma we will be writing $c_i$ instead of $c_{i, n}$.

\bigskip 

 Recall the general setup of Nakajima's Chapter 9 (or rather its
adaptation to our case). We choose and fix two smooth transversal 
curves $C$ and $C'$ in one very ample linear  system on $S$. 
Take such an $s$  that $M^G(r, s)$ is non-empty  and fix a 
locally free sheaf ${\cal E}$ correspoding to  a point in $M^G(r, s)$.

We   will  compute numbers   $c_i$ using  the  formalizm of generating
functions   and a connection with   symmetric  functions discovered by
Nakajima.  The actual idea of using  embedded curves is originally due
to Grojnowski.

\bigskip

To that end, consider the subvariety $L^{k, s} \subset M^G(r, n)$  
formed by all  sheaves $\F$ such that

 (a) $\F^{**} \simeq {\cal E}$;

 (b) The quotient $A_{\F}:= {\cal E}/\F$ is  an Artin sheaf of length 
$k = n-s$ supported at finitely many points of $C$;

\medskip 

{\bf Definition.}  For each partition $\mu= (\mu_1 \geq \ldots \geq \mu_l 
> 0)$ of $k$ let $(L^{\mu, s})^{\circ}$ be the set of all sheaves $\F$ such 
that

(a) $\F^{**} \simeq {\cal E}$;

(b) the quotient $A_{\F} = {\cal E}/\F$ is supported at pairwise distinct
points $x_j \in C$, $1 \leq j \leq l$ with multiplicities $\mu_j$.

Denote also by $L^{\mu, s}$ the closure of $(L^{\mu, s})^{\circ}$.

\bigskip

{\bf Definition.} Similarly, we define $\widehat{L}^{\mu, s}$ to be
 the closure of the subset of all sheaves $\F$ satisfying:

(a) $\F \subset \F'$ for some $\F' \in M^G(r, s)$;

(b) the quotient $\F'/\F$ is supported at pairwise distinct
points $x_j \in C$, $1 \leq j \leq l$ with multiplicities $\mu_j$.

\bigskip

\begin{prop} The closed subvarieties $L^{\mu, s}$ are the irreducible
components of $L^{k, s}$. Each component $L^{\mu, s}$ is of dimension 
$rk$. Similarly, $\widehat{L}^{\mu, s}$ is the set of all irreducible
components of $\widehat{L}^{k, s}$ and all  $\widehat{L}^{\mu, s}$ are
of the same dimension.
\end{prop}
{\it Proof.} If we fix all the points $x_j$ from the definition of 
$L^{\mu, s}$ then by Theorem 2.2 the variety of all sheaves with the 
required  condition has dimension $\sum_{i=1}^l (r\mu_i - 1) = rk -l$. 
Since the points $x_j$, $j = 1, \ldots, l$ are allowed to move in an 
$l$-dimensional family, we obtain $\dim L^{\mu, s} = rk$. By the 
irreducibility part of Theorem 2.2 we conclude that each $L^{\mu, s}$ 
is irreducible.

It follows from the definition of $L^{k, s}$ that it is a union
of the closed subsets $L^{\mu, s}$. Hence, $L^{\mu, s}$ are the 
irreducible components of $L^{k, s}$. \quad $\Box$

\bigskip 

Note that $P^C_{-i}$ maps classes of (cohomological) degree $b$ on
$M^G(r, n)$ to classes of degree $ri + b$ on $M^G(r, n+i)$. In
particular, the subspace $\bigoplus H^{rn}(M^G(r, n))$ is preserved by
$P^C_{-i}$. Moreover, the next theorem shows that in fact even the 
subspace generated by the classes of $L^{\mu, s}$, is invariant under 
$P^C_{-i}$ ($i > 0$).

\begin{theorem} (cf. [Na, Theorem  9.14])
For $i > 0$, we have  $P^C_{-i}[L^{\mu, s}] = \sum_\lambda a_{\lambda
  \mu} [L^{\lambda, s}]$,
where the summation is over partitions $\lambda$ of $|\mu| + i$
which are obtained as follows: 

(a) add $i$ to a term in $\mu$, say $\mu_j$
(possibly 0), and then 

(b) arrange it in descending order. 

\noindent The coefficient
$a_{\lambda\mu}$ is $\#\{l\;|\; \lambda_l = \mu_j + i\}$.

A similar statement is true for the classes of $\widehat{L}^{\mu, s}$.
\end{theorem}

\bigskip

This theorem will be proved in Section 5.

\hspace{2cm}

Let  $[vac]$ be the cohomology  class of the  point $\E \in M^G(r, s)$
and $[Vac]$ be  the fundamental  class of  $M^G(r,  s)$.  Theorem  4.6
allows  us  to establish  a  connection with  the theory  of symmetric
functions as follows. Let $L$ (resp. $\widehat{L}$) be the subspace in
$\bigoplus_n H^*(M^G(r, n))$ generated by  the classes of $L^{\mu, s}$
and $[vac]$ (resp. $\widehat{L}^{\mu,   s}$ and $[Vac]$). We  define a
${\Bbb C}$-linear isomorphism from  $L$ (resp. $\widehat{L}$) onto the
space $\Lambda$ of symmetric   functions in infinitely  many variables
(cf.  [Na,  9.1],  [Mc2]).  This    isomorphism sends  $[vac]$ (resp.  
$[Vac]$)   to   $1  \in   \Lambda$    and   $[L^{\mu, s}]$    (resp.   
$[\widehat{L}^{\mu, s}]$) to the orbit sum function $m_{\lambda}$ (cf.
{\it  loc. cit.}).

Theorem 4.6 means that the  operator  $P^C_{-i}$ corresponds under  the
isomorphism above to multiplication by the $i$-th power sum (or Newton
function) $p_i \in \Lambda$ ({\it  loc. cit.}).  This provides us with
several identities  between cohomology  classes coming  from classical
identities between symmetric functions.

  Note that for $\mu=(1^k)$  the corresponding variety $L^{\mu, s}$ is
isomorphic to the global Quot scheme $Quot_C^k(\E)$ on the curve $C$. 
This scheme parametrizes quotient sheaves ${\cal E}|_C \to A$ on $C$, 
where $A$ is of length $k$. Every such quotient on $C$ defines a sheaf 
on $S$, namely the kernel of the composition ${\cal E} \to {\cal E}|_C 
\to A$. The class of $Quot^k_C(\E)$ corresponds to the $k$-th elementary
symmetric function $e_k = m_{(1^k)}$.

 Similarly, we will slightly abuse notation by writing $Quot^k_C(s)$
instead of $\widehat{L}^{(1^k), s}$. This variety is a birational
image of the family of schemes $Quot_C^k(\F')$ parametrized by 
points $\F'$ of $M^G(r, s)$.

\bigskip 

Repeating the arguments of Nakajima, (see [Na, formula 9.16]) we
obtain the formulas
\begin{equation}
\label{exp}
\sum_{n=0}^{\infty} z^n [Quot_C^n(\E)] = exp \Big( \sum_{i = 1}^{\infty} 
{z^i P^C_{i} \over (-1)^{i-1} i } \Big) \cdot [vac].
\end{equation}
$$
\sum_{n=0}^{\infty} z^n [Quot_C^n(s)] = exp \Big( \sum_{i = 1}^{\infty} {
  z^i P^C_{i} \over (-1)^{i-1} i } \Big) \cdot [Vac].
$$

These    formulas  arise  from the  classical   identity between Newton
symmetric functions and  elementary  symmetric functions (cf.  [Mc2]). 
They are   the first ingredient in  our  computation of   $c_{i}$. The
second ingredient is provided by

\begin{theorem}(cf. [Na, Exercise 9.23]) The following relation holds:
$$
\sum_{n=0}^{\infty} z^{2n} \langle Quot^n_C(\E), Quot^n_{C'}(s) \rangle =
(1 - (-1)^r z^2)^{r \langle C, C' \rangle}.
$$
\end{theorem}

The proof of this theorem will be given in Section 6 for the case when
the intersection $C \bigcap C'$ is transversal (which is sufficient
for our purposes).

\bigskip

{\it End of Proof of Theorem 4.1:} Following Nakajima, we
 introduce the following  notations:
$$
C_-(z) = \sum_{i=1}^{\infty} {P^C_{-i} z^i \over (-1)^{i-1} i};
\quad \quad \quad 
C_+(z) = \sum_{i=1}^{\infty} {P^C_{i} z^i \over (-1)^{i-1} i}.
$$
Note that $C_-(z)$ is adjoint to $C_+(z)$ with respect to the 
intersection form on $\bigoplus_n H^*(M^G(r, n))$ (the cohomology groups
for different $n$ are orthogonal to each other). We extend this
form to power series in $z$ by $z$-linearity.
By Theorem 4.7 and (\ref{exp}) one has:
$$
(1  - (-1)^r z^2)^{r \langle C,   C' \rangle} = \sum_{n=0}^{\infty}
z^{2n} \langle Quot^n_C(\E), Quot^n_{C'}(s) \rangle =
$$
$$ =
\langle  \sum_{n=0}^{\infty} z^{n}
Quot^n_C(\E),\; \sum_{n=0}^{\infty} z^{n} Quot^n_{C'}(s) \rangle = 
$$
$$
=
 \langle exp(C_-(z)) \cdot [vac],\; exp(C'_-(z)) \cdot [Vac] \rangle =
$$
$$ =
 \langle exp(C'_+(z))\;  exp (C_-(z)) \cdot [vac], [Vac] \rangle 
$$
where the last equation follows from adjointness.

 For any pair of operators $A$ and $B$, one has an identity 
$$
exp(-A)\; B\; exp (A) = exp(-ad\; A) \big(B \big) = 1 - (ad\; A)(B) + 
{(ad\; A)^2 \over 2!}(B) - \ldots
$$
Therefore
$$
\langle exp(C_+(z))\; exp (C_-(z)) \cdot [vac],\; [Vac] \rangle = \hspace{5cm}
$$
\hfill $$
\hspace{2cm}= 
\langle exp (C_-(z)) \Big[ exp( - ad\; C_-(z)) \Big(exp(C'_+(z)) \Big) \Big]
\cdot [vac], [Vac] \rangle
$$ 

\noindent
An explicit computation shows that 
$$
\big[ C_-(z),\; exp(C_+(z)) \big] = - \Big(\sum_{n=1}^{\infty} {c_n \over n^2}
\langle C, C' \rangle z^{2n} \Big)\; exp(C_+(z))
$$
Denote the  expression  $\displaystyle  \sum_{n=1}^{\infty} { c_n
  \over n^2} \langle C, C' \rangle z^{2n} $  by $\Phi(z)$. Then by the
previous formula:
$$
exp( -     ad\;   C_-(z)) \Big(exp(C_+(z))   \Big)=    exp(\Phi(z))
exp(C_+(z)).
$$
 Collecting the results of computations, we obtain:
$$
(1 - (-1)^r z^2)^{r \langle C, C' \rangle} = \hspace{6cm}
$$
$$=\langle exp (C_-(z)) 
\Big[ exp( - ad\; C_-(z)) \big(exp(C_+(z)) \big) \Big]
\cdot [vac], [Vac] \rangle=
$$
$$
= exp(\Phi(z)) \langle exp (C_-(z))\; exp(C_+(z)) \cdot [vac], [Vac] \rangle
= 
$$
$$
= exp(\Phi(z)) \langle exp (C_-(z)) \cdot [vac], [Vac] \rangle
= exp(\Phi(z))
$$
The last equality holds since all $P^C_{-i}$ involved in the
definition of $C_+(z)$, map  $[vac]$ to the 
orthogonal complement of $[Vac]$. By definition of $\Phi(z)$
we have
$$
\sum_{n=1}^{\infty} {c_n \over n^2} \langle C, C' \rangle  z^{2n} 
= r \langle C, C' \rangle  log(1 - (-1)^r z^2).
$$
Hence $c_n = (-1)^{rn-1} rn$ which completes the proof of Theorem 4.1 $\Box$

\section{A transversality result.}

  In $r=1$ case Theorem 4.6 has a simple proof (cf. [Na]) using local
coordinates on $S$. While this result is still true in higher ranks
the proof of it requires a more detailed analysis of the geometry
of the moduli space to be provided below.

  First consider the situation at the set-theoretic level.
It follows for the definitions that if $(\F_1, \F_2)  \in P^C_{-i}$ 
then $\F_1^{**} \simeq  \F_2^{**}$ and $Supp(\F_1/\F_2) = \{x\} \in C$. 
In particular, if
$\F_1 \in L^{\mu, s}$ then $\F_2$ necessarily belongs to one of the
$L^{\lambda, s}$ where the partition $\lambda$ is as described in the
theorem. Hence $P^C_{-i}[L^{\mu, s}]$ is a linear combination of 
$[L^{\lambda, s}]$ with integral coefficients.

\bigskip

Our next step is to show that for generic $\F_2 \in L^{\lambda, s}$ 
the intersection $$p_1^{-1}(L^{\mu, s})
\bigcap P^C_{-i}  \bigcap p_2^{-1}(\F_2) \subset M^G(r, s + |\mu|)
\times  M^G(r, s+ |\lambda|)$$ is finite and consists of 
exactly $a_{\lambda \mu}$ points.

Consider the quotient $A_{\F_2} = {\cal E}/\F_2$. If $\lambda=
(\lambda_1 \geq \ldots \lambda_p > 0)$ then $A_{\F_2}$ is supported
at some points $x_j \in C$, $1\leq j \leq p$ with multiplicities 
$\lambda_j$. Therefore one can write $A_{\F_2}$ as a direct sum 
$\bigoplus (A_{\F_2})_{x_j}$.

We claim that for generic $\F_2$ there exist local 
coordinates $(\zeta_j, \xi_j)$ at $x_j$ such that $(A_{\F_2})_{x_j} \simeq
{\Bbb C}[\zeta_j]/(\zeta_j)^{\lambda_j}$.

In fact, choose a trivialization ${\cal E} \simeq {\cal O}^{\oplus r}$
in  the neighbourhood of  $x_j$  and fix   some system of  coordinates
$(\zeta'_j, \xi'_j)$ centered at $x_j$. Consider such quotients ${\cal
  O}^{\oplus r} \to  A$ of length $\lambda_j$  at $x_j$ that the first
component ${\cal O}  \to A$ is  surjective and its kernel is generated
by  $m_{x_j}^{\lambda_j}$ and    some  element   $\xi_j = \xi'_j     +
\sum_{k=1}^{\lambda_j - 1}  a_k   (\zeta_j')^k$, where  $a_1,  \ldots,
a_{\lambda_j -1}$ are arbitrary complex numbers.  The other components
${\cal O}^{\oplus (r-1)} \to A$ can be chosen  arbitrarily. If we take
$\zeta_j = \zeta_j'$  then  $A$ has the   required form  in  the local
coordinates  $(\zeta_j,  \xi_j)$.  Thus    we  obtain  a  $(r\lambda_j
-1)$-dimensional  family of non-isomorphic quotients.  Here $\lambda_j
-1$  parameters  come   from  the choice    of $\xi_j$ and   the other
$(r-1)\lambda_j$ parameters come from  the  choice of the map   ${\cal
  O}^{\oplus (r-1)} \to  A$. By Theorem 2.2  this family forms a dense
subset in $Quot(r, \lambda_j)$ and this implies our assertion.

Hence  we can assume  that ${\cal  E}/\F_2$  satisfies the  conditions
described above. Then a  choice of $\F_1$ such  that $(\F_1, \F_2) \in
P^C_{-i}$  amounts to choosing  a $(\lambda_j -i)$-dimensional quotient
$(A_{\F_1})_{x_j}$   of  $(A_{\F_2})_{x_j}   \simeq {\Bbb C}[\zeta_j]/
(\zeta_j)^{\lambda_j}$. If  $x_j$  is fixed,  this  can  be  done in a
unique way.  Therefore the number  of points in  $p_1^{-1}(L^{\mu, s})
\bigcap P^C_{-i}$ over a generic $\F_2$ is equal  to the number of ways
to  subtract  $i$ from one  of  the parts of  partition $\lambda$, and
obtain partition $\mu$.  This number is exactly $a_{\lambda \mu}$.

\bigskip 

 To finish the proof we need to show that $p^{-1}(L^{\mu, s})$
  intersects  $P^C_{-i}$ transversly.  To that end, we
prove the following lemma

\begin{lemma} Let $(\F_1, \F_2) \in P^C_{-i} \bigcap M^G(s+ |\mu|) \times 
  M^G(s + |\lambda|)$ be as smooth point of  $P^C_{-i}$ and assume that
  $\F_1 \in L^{\mu, s}$ (resp.  $\F_2 \in L^{\lambda, s}$) is  generic
  in the sense described above. Then the intersection
$$
  W =   T_{(\F_1, \F_2)}\;p_1^{-1}(L^{\mu,  s})  \bigcap  T_{(\F_1,
    \F_2)}\;P^C_{-i}
$$
of tangent spaces to $p_1^{-1}(L^{\mu, s})$ and $P^C_{-i}$ 
projects isomorphically under $p_2$ onto the tangent space 
$T_{\F_2}(L^{\lambda, s})$.
\end{lemma}
{\it Proof.} It suffices to prove
that if $(v, w) \in W$ then $w \in T_{\F_2}(L^{\lambda, s})$ 
and $v$ is uniquely defined by $w$. Then the dimension count shows that in 
fact the map $dp_2: W \to T_{\F_2}(L^{\lambda, s})$ is an isomorphism.
Our proof consists of several steps.

\bigskip

{\it Step 1.} First we compute the tangent spaces $T_{\F_1}(L^{\mu, s})$ 
and $T_{\F_2}(L^{\lambda, s})$.

To that end, recall (cf. [Ar]) that the tangent space $T_{\F_1}(M^G(r, s+ |\mu|))$
is isomorphic to the kernel of the natural trace map $tr^1: Ext^1_S (\F_1, \F_1) \to 
H^1(S, {\cal O})$. Infinitesimal deformations of $\F_1$ with fixed ${\cal E}
= \F_1^{**}$ correspond to the subspace $Hom_S (\F_1, {\cal E}/\F_1)
\subset Ext^1_S (\F_1, \F_1)$. The embedding of this subspace  
is just the  boundary map which comes from applying  
$Hom(\F_1, \;\cdot\;)$ to a short
exact sequence $0 \to \F_1 \to {\cal E} \to {\cal E}/\F_1 \to 0$.

  Since $\displaystyle {\cal E}/\F_1 = A_{\F_1} = \bigoplus_{x_j 
\in Supp(A_{\F_1})} (A_{\F_1})_{x_j}$, one
has a direct sum decomposition 
$$
Hom_S (\F_1, A_{\F_1}) = \bigoplus_{x_j \in Supp(A_{\F_1})}
 Hom_S (\F_1, (A_{\F_1})_{x_j}).
$$
Recall that $\mu_j = mult_{x_j} (A_{\F_1})$. Suppose that $\F_1$ is 
generic in the sense that it has a local description given before this
lemma. It is easy to show using this local description that 
$T_{\F_1}(L^{\mu, s})$ corresponds precisely to the subspace
$$
\bigoplus_{x_j \in Supp(A_{\F_1})} Hom_S (\F_1/ {\cal E}(-\mu_j C), 
(A_{\F_1})_{x_j}) \subset Hom_S (\F_1, A_{\F_1}).
$$
Similar computation applies to $\F_2 \in L^{\lambda, s}$.

\bigskip 

{\it Step 2.} 
Recall that any tangent vector $v \in T_{\F_1}(M^G(r, s+ | \mu|)) \subset
Ext^1_S(\F_1, \F_1)$ corresponds
to an extension $0 \to \F_1 \to {\cal G}_1 \to \F_1 \to 0$ (i.e.
deformation with base ${\Bbb C}[\epsilon]/\epsilon^2$). 

The conditions $Supp\; (\F_1/\F_2) = x \in C$ and
$length(\F_1/\F_2) = i$ imply that $\F_1(-iC) \subset \F_2$. Moreover,
 the following diagram 

\begin{equation}
\label{diag}
\begin{CD}
\F_1(-iC) @>>> \F_2 \\
@{|} @VVV \\
\F_1(-iC) @>a>> \F_1
\end{CD}
\end{equation}
commutes.  Here the map $a$ is the muplitplication by  $i$-th power
of the local equation for $C$.

 The diagram (\ref{diag}) should be preserved under infinitesimal 
deformation of the point $(\F_1, \F_2)$ within $P_C[-i]$. This
means that

(a) There exists a commutative diagram:
$$
\begin{CD}
0   @>>>   \F_2 @>>>  {\cal G}_2 @>>>   \F_2 @>>>  0 \\
  &  &      @VVV        @VVV           @VVV     & \\
0   @>>>   \F_1 @>>>  {\cal G}_1 @>>>   \F_1 @>>>  0 
\end{CD}
$$
It is a standard fact of homological algebra that this is equivalent to 
requiring that the images of $v$ and $w$ under the natural maps 
$$Ext^1_S (\F_1, \F_1)
\to Ext^1_S (\F_2, \F_1) \quad \mbox{ and }\quad Ext^1_S (\F_2, \F_2) 
\to Ext^1_s 
(\F_2, \F_1)$$
coincide.

(b) There exists a similar diagram of extensions corresponding
to  the embedding $\F_1(-iC) \subset \F_2$. This is can be
expressed in terms of Ext groups in a similar way.

(c) The diagram of middile terms 
$$
\begin{CD}
{\cal G}_1(-iC) @>>> {\cal G}_2 \\
@{|} @VVV \\
{\cal G}_1(-iC) @>a>> {\cal G}_1 
\end{CD}
$$
commutes. This condition can be expressed as  vanishing of some
homomorphism from the quotient copy of $\F_1(-iC)$ to the subsheaf copy 
of $\F_1$ (i.e. the embedding of sheaves $a$ should not be deformed). 
We will not write this down explicitly as we need this condition only 
in a special case (see below).

\bigskip

{\it Step 3.} We will prove that if $(v, w) \in U$ and $v \in T_{\F_1}
(L^{\mu, s})$ then $w  \in T_{\F_2}(L^{\lambda, s})$ and $w=0$ implies
$v=0$.    The condition $w   \in T_{\F_2}(L^{\lambda, s})$ will follow
from (a) above while the   implication $(w=0) \Rightarrow (v=0)$ is  a
consequence of (c).

  To use (a), consider the diagram:
\begin{equation}
\label{def-diag}
\begin{CD}
0   @>>>   \F_1 @>>>  {\cal E} @>>>   A_{\F_1} @>>>  0 \\
  &  &      @AAA        @AAA           @AAA     & \\
0   @>>>   \F_2 @>>>  {\cal E} @>>>   A_{\F_2} @>>>  0 
\end{CD}
\end{equation}
and the induced commutative digram of Ext-groups:
$$
\begin{CD}
  Hom_S(\F_1, A_{\F_1}) @>>>  Ext^1_S(\F_1, \F_1) @>>> Ext^1_S (\F_1, {\cal E})
 \\
   @VVV        @VVV           @VVV      \\
  Hom_S(\F_2, A_{\F_1}) @>>>  Ext^1_S(\F_2, \F_1) @>>> Ext^1_S (\F_2, {\cal E})
 \\
   @AAA        @AAA           @{|}      \\
  Hom_S(\F_2, A_{\F_2}) @>>>  Ext^1_S(\F_2, \F_2) @>>> Ext^1_S (\F_2, {\cal E})
  \\
\end{CD}
$$
\medskip
\noindent
Here the  lower two   rows   are obtained  by  applying   $Hom_S(\F_2,
\;\cdot\;)$ to (\ref{def-diag}) and  the upper row comes brom applying
$Hom_S(\F_2,  \;\cdot\;)$ to the upper row  of (\ref{def-diag}).  Note
that by stability  $Hom_S(\E, \E) =  Hom_S (\F_i,  \F_i) = Hom_S(\F_i,
\E) = \C$ and hence the first arrow in  each row is injective.  Simple
diagram chase shows that if
$$
v   \in   \bigoplus_{x_j \in  Supp(A_{\F_1})}   Hom_S  (\F_1/ {\cal
  E}(-\mu_j    C), (A_{\F_1})_{x_j})  \subset  Hom_S (\F_1,  A_{\F_1})
$$
 then by the condition (b) we have
$$
 w \in   \bigoplus_{x_j   \in Supp(A_{\F_2})} Hom_S   (\F_2/  {\cal
   E}(-\lambda_j C), (A_{\F_2})_{x_j})  \subset Hom_S (\F_2, A_{\F_2}).
$$

\bigskip 

To prove injectivity, recall that  if $v \in Hom_S(\F_1, A_{\F_1})
\subset  Ext^1(\F_1, \F_1)$ then the  extension $$0 \to  \F_1 \to {\cal
  G}_1 \to \F_1 \to 0$$  can  be recovered as a   kernel of the map  of
extensions
$$
\begin{CD}
0 @>>> {\cal E} @>>> {\cal E} \oplus \F_1 @>>> \F_1 @>>> 0 \\
& & @VbVV @V(b \oplus v)VV @VVV \\
0 @>>> A_{\F_1} @{=} A_{\F_1} @>>> 0 @>>> 0 
\end{CD}
$$
where $b$ is the cokernel of the natural embedding $\F_1 \subset {\cal E}$.
Similar statement is true for $w$ and the corresponding
extension $0 \to \F_2 \to {\cal G}_2 \to \F_2 \to 0$.
One shows that  in our case the condition (c) translates into 
commutativity of the diagram
$$
\begin{CD}
\F_1(-iC) @>v>> A_{\F_1}(-iC) \\
@VVV @VVV \\
\F_2 @>w>> A_{\F_2}
\end{CD}
$$
Since both vertical arrows in this diagram are injective the condition $w=0$
implies $v =0$. 

\bigskip

 This completes the proof of the lemma and hence the proof of 
Theorem 4.6.
\quad $\Box$

\section{Computation of the intersection number.}

  The main result of this section is a computation of the 
intersection number of the fundamental classes of cycles 
$Quot^n_C(\E)$ and $Quot^n_{C'}(s)$  in $M^G(r, n+s)$.
In the case of Hilbert schemes this intersection is trasversal 
and the intersection number can be found by a simple set-theoretic
argument (recall that we assume that $C$ and $C'$ are transversal).
However, for higher ranks the set-theoretic intersection is not
transversal any more and to compute the intersection index we have to
apply the excess intersection formula (cf. [Fu]).

  By definition $Quot^n_C(\E) \subset M^G(r, n+s)$ parametrizes all 
short exact sequences $0 \to \E_1 \to \E|_C \to A \to 0$ of sheaves on
$C$, where $A$ is a length $n$ Artin sheaf. Any such sequence defines 
a sheaf $\F \in M^G(r, n+s)$ on $S$, namely the kernel of the composition
$\E \to \E|_C \to A$. Note that $Quot^n_C(\E)$ is {\it smooth} since
the tangent space to it at a point represented by $0 \to \E_1 \to
\E|_C \to A \to 0$ is $Hom_{\Oo_C}(\E_1, A)$. Since $\E_1$ is a rank r
vector bundle on $C$ this space is of constant dimension $rn$. 

Similarly,  $Quot^n_{C'}(s)   \subset M^G(r,  n+s)$   has open  subset
$(Quot^n_{C'}(s))^{\circ}$ isomorphic   to  a fiber  bundle over   the
subset  $(M^G(r, s))^{\circ}$   of  locally free   sheaves, with fiber
$Quot^n_{\C'}(\E')$   over  $\E'  \in  (M^G(r,   s))^{\circ}$.   Hence
$(Quot^n_{C'}(s))^{\circ}$ is smooth.

\bigskip

Note that $Quot^n_{C'}(s) \setminus (Quot^n_{C'}(s))^{\circ}$ does not
intersect $Quot_C^n(\E)$ since all $\F \in Quot^n_C(\E) \subset M^G(r,
n+s)$ satisfy  $length(\F^{**}/\F) =      n$  and  for    sheaves     in
$Quot^n_{C'}(s) \setminus  (Quot^n_{C'}(s))^{\circ}$ this length is at
least $(n+1)$.

Assume for the  sake of simplicity  that $S =  {\Bbb P}^2$ and $C, C'$
are two distinct lines intersecting at a point $x  \in S$. The general
case follows from our proof by a simple combinatorial argument.

\bigskip

Let $\F \in Quot^n_{C'}(s) \bigcap Quot^n_C(\E)$. Then $\F^{**} = \E$,
$\F$ is a kernel of $\E \to \E|_C \to A$ and also a kernel of 
$\E \to \E|_{C'} \to A'$. Hence $\F$ contains $\E(-C) + \E(-C') =
\E \otimes {\cal J}_x$ as a subsheaf and $A \simeq A' \simeq
(\C_x)^n$, where ${\cal J}_x$ is the ideal sheaf of $x \in S$ and $\C_x 
\simeq \Oo_S/{\cal J}_x$. This means that every sheaf $F$ in our
intersection can be obtained as a kernel  of $\E \to \E_x \to A$ where
$A$ is an $n$-dimensional quotient of an $r$-dimensional vector space
$\E_x$. Hence $Quot^n_{C'}(s) \bigcap Quot^n_C(\E)$ is nothing but the
Grassman variety $Gr(\E_x, n)$ of all $n$-dimensional quotients of 
$\E_x$. In particular, $n \leq r$.

\bigskip

 Since $Quot^n_{C'}(s)$ and $Quot^n_C(\E)$ are of complementary 
dimensions in the ambient $M^G(r, n+s)$ we see that for $r \geq 2$ their 
set-theoretic intersection has abnormally high dimension. However,
it is true that $Quot^n_{C'}(s)$ and $Quot^n_C(\E)$  are smooth
along the points of $Gr(\E_x, n)$ and, moreover, for any $\F \in
Gr(\E_x, n)$ 
$$
T_{\F}\; Gr(\E_x, n) =  T_{\F}\; Quot^n_{C'}(s) \bigcap T_{\F}\; 
Quot^n_C(\E)
$$
(this can be checked using methods of the previous section).

Therefore, we can apply excess intersection formula (cf. [Fu, Example
6.1.7]): let $V$ be a rank $(r-n)n$ vector bundle on $Gr(\E_x, n)$
arising form the following exact sequence: 
\begin{equation}
\label{excess}
0 \to T_{Gr(r, n)} \to T_{Quot_{C}^n(\E)} \oplus T_{Quot_{C'}^n(s)} 
\to T_{M^G(r, n+s)} \to V \to 0. 
\end{equation}
Then the intersection number $[Quot_{C}^n(\E)] \cdot
[Quot_{C'}^n(s)]$ in $M^G(r, n+s)$ is equal to the top Chern class
$c_{(r-n)n} (V)$.

\bigskip

To compute this Chern  class we will need  to consider certain sheaves
on $M^G(r, n+s)   \times S$. For any  pair  of sheaves  ${\cal  G}_1$,
${\cal G}_2$ on $M^G(r, n+s) \times S$ let $\E xt^i_{p_1} ({\cal G}_1,
{\cal  G}_2)$ be  the relative  Ext-sheaf on $M^G(r, n+s)$ with   
respect to the  first
projection $p_1: M^G(r, n+s)  \times S \to M^G(r,  n+s)$. We will only
deal with the cases  when for all  $x \in M^G(r,  n+s)$ the global Ext
group  $Ext^i(p_1^{-1}(x);  {\cal  G}_1, {\cal   G}_2)$   on the fiber
$p_1^{-1}(x)$ is of  constant dimension. Hence  by [Kl] $\E xt^i_{p_1}
({\cal G}_1, {\cal G}_2)$ is  a vector bundle  on $M^G(r, n+s)$ of the
same rank.  Similar remarks apply  to  any closed subspace of  $M^G(r,
n+s)$ (abusing notation we will denote a sheaf on $M^G(r, n+s)$ and its
restriction to a closed subspace by the same letter).

\bigskip

Consider a larger moduli space $\widetilde{M}^G(r, n+s)$ of $H$-stable
sheaves   with $c_1(\F)  =  c_1(L)$   and  $c_2(\F) =   n+s$ (i.e. the
determinant is not fixed but is allowed to differ from $L$ by any line
bundle in $Pic^0(S)$). Since $gcd(r, c_1(H)  \cdot c_1(L)) = 1$, there
exists a universal sheaf ${\cal G}$ on $\widetilde{M}^G(r, n+s) \times
S$  (cf. [HL]). It  is known (cf. [Ma])  that  the tangent bundle $T\;
\widetilde{M}^G(r,  n+s)$ is  isomorphic  to $\E xt^1_{p_1}({\cal  G},
{\cal G})$.  Since  $M^G(r, n+s)$  is a fiber   of a smooth  map $det:
\widetilde{M}^G(r, n+s) \to Pic^0(S)$,  the tangent bundle  to $M^G(r,
n+s)$  has  the same Chern classes  as  $\E xt^1_{p_1}({\cal G}, {\cal
  G})$.

\bigskip

To compute  the full  Chern class  of  $\E xt^1_{p_1}({\cal  G}, {\cal
  G})|_{Gr(\E_x,  n)}$  we use a   short exact sequence  of sheaves on
$Gr(\E_x, n) \times S$:
\begin{equation}
\label{seq}
0 \to {\cal G} \to p_2^*\E \to {\cal A} \to 0,
\end{equation}
where $p_1$,  $p_2$ are the projections,  ${\cal A}  = p_1^* Q \otimes
p_2^*  \C_x$   and $Q$  is   the universal   quotient  bundle  on  the
Grassmanian.

\bigskip 

First of  all,  applying  $R{\cal   H}om_{p_1}(\cdot,  {\cal  G})$  to
(\ref{seq}) one   obtains  a long  exact  sequence  of  sheaves on the
Grassmanian:
$$
0  \to  {\cal H}om_{p_1}({\cal  G},  {\cal  G}) \to  \E xt^1_{p_1}
({\cal A}, {\cal  G}) \to \E xt^1_{p_1} (p_2^* \E,
{\cal G})   \to \hspace{3cm} 
$$
$$
\hspace{6cm} \to 
\E xt^1_{p_1}({\cal   G},  {\cal   G})  \to \E
xt^2_{p_1}({\cal A}, {\cal G}) \to 0
$$

Note that ${\cal H}om_{p_1}(p_2^* \E,  {\cal  G})$ vanishes since its
fiber over $\F \in Gr(\E_x, n)$ is $Hom_{\Oo_S}(\E, \F)$ which is zero
($\E$ and $\F$ are stable and $c_2(\F) > c_2 (\E)$). Moreover, the
fiber of ${\cal H}om_{p_1}({\cal  G},  {\cal  G})$ over $\F$ is 
$Hom_{\Oo_S}(\F, \F)$. Since stable sheaves have only scalar
automorphisms, ${\cal H}om_{p_1}({\cal  G},  {\cal  G})$ is a trivial
line bundle. One can also check that either both $ \E xt^2_{p_1}( p_2^* \E,
{\cal G})$ and   $\E xt^2_{p_1}({\cal G}, {\cal G})$  are zero (when
$c_1(K) \cdot c_1(H) < 0$) or the map between them is an isomorphism (when $K
\simeq \Oo$).

  Thus, the full Chern class of $T_{M^G(r, n+s)}$ is equal to:
$$
  c(T_{M^G(r, n+s)})  = c(\E xt^1_{p_1}({\cal  G},  {\cal G})) =  {
    c(\E  xt^2_{p_1}({\cal A}, {\cal G}))   c( \E xt^1_{p_1}(p_2^* \E,
    {\cal G})) \over c(\E xt^1_{p_1} ({\cal A}, {\cal G}))}
$$

\vspace{0.7cm}

\noindent
As a second step we apply $R{\cal H}om_{p_1}(p_2^*\E, \cdot)$ to
(\ref{seq}) and get

$$
0 \to {\cal H}om_{p_1} (p_2^* \E, p_2^* \E) \to {\cal H} om_{p_1}
(p_2^* \E, {\cal A}) \to \hspace{4cm}
$$
$$
\hspace{5cm} \to  \E xt^1_{p_1} (p_2^* \E, {\cal G}) \to
\E xt^1_{p_1} (p_2^* \E, p_2^* \E) \to 0
$$

\noindent
Here, the sheaf ${\cal H} om_{p_1}(p_2^* \E, {\cal G})$ vanishes as
before
and $\E xt^1(p_2^* \E, {\cal A})$ is zero since its 
fiber over a point $\F$ is $Ext^1_{\Oo_S}(\E, (\C_x)^n) = H^1(S,
\E^* \otimes (\C_x)^n) = 0$.  The bundles ${\cal H}om_{p_1} (p_2^* \E,
p_2^* \E)$ and $\E xt^1_{p_1} (p_2^* \E, p_2^* \E)$ are trivial on
$Gr(\E_x, n)$ hence

$$
c( \E xt^1_{p_1}(p_2^* \E, {\cal G})) = c({\cal H} om_{p_1} (p_2^* \E,
{\cal A})) = c(Q \otimes \C^r) = (c(Q))^r.
$$

\vspace{0.7cm}

\noindent
Finally, we apply $R{\cal H}om_{p_1}({\cal A}, \cdot)$
to (\ref{seq}) and get

$$
0 \to {\cal H}om_{p_1} ({\cal A}, {\cal A}) \to Ext^1_{p_1} 
({\cal A}, {\cal G}) \to 0 \to Ext^1_{p_1} ({\cal A}, {\cal A}) \to 
$$
$$
\to  Ext^2_{p_1} ({\cal A}, {\cal G}) \to  Ext^2_{p_1} ({\cal A},
p_2^* \E) \to Ext^2_{p_1} ({\cal A}, {\cal A}) \to 0
$$

\noindent
where $\E xt^1_{p_1}({\cal A}, p_2^* \E)$ vanishes since it is dual
to $\E xt^1_{p_1}(p_2^* \E, {\cal A} \otimes p_2^* K_S)$ (cf. [Kl]).
Since $Hom_{\Oo_S}(\C_x, \C_x) = Ext^2_{\Oo_S}(\C_x, \C_x) = \C$ and
$Ext^1_{\Oo_S}(\C_x, \C_x) = \C^2$, we have
$$
{\cal H}om_{p_1} ({\cal A}, {\cal A}) \simeq Ext^2_{p_1} 
({\cal A}, {\cal A}) \simeq Q \otimes Q^*, \quad 
Ext^1_{p_1} ({\cal A}, {\cal A}) = (Q \otimes Q^*)^{\oplus 2}.
$$

\vspace{0.7cm}
\noindent
To compute $\E xt^2_{p_1} ({\cal A}, p_2^* \E)$ we again apply a
relative version of Serre's duality ({\it loc. cit.}) which gives

$$
\E xt^2_{p_1} ({\cal A}, p_2^* \E) \simeq (Q^*)^{\oplus r}.
$$

Summing up the results of our computation, we obtain

$$
c(T_{M^G(r, n+s)}) = {(c(Q^*))^r (c(Q))^r (c(Q \otimes Q^*))^2 \over
c(Q \otimes Q^*) c(Q \otimes Q^*)} = (c(Q) c(Q^*))^r.
$$

A similar approach can be used with $T_{Quot^n_C(\E)}$ and
$T_{Quot^n_{C'}(s)}$ (one has to consider sheaves on $Gr(\E_x, n)
  \times C$ and $Gr(\E_x, n) \times C'$). One checks that 
$$
c(T_{Quot^n_C(\E)}) = c(T_{Quot^n_{C'}}(s)) = (c(Q))^r.
$$

\bigskip 

Let $S$ be the universal subbundle on $Gr(\E_x, n)$. Using the exact sequence
(\ref{excess}) and  the following short exact sequence: 
$$
0 \to Q \otimes Q^* \to Q^{\oplus r} \to S^* \otimes Q \to 0
$$
(recall that $T_{Gr(\E_x, n)} \simeq S^* \otimes Q$), we get
$$
c(V) = {c(T_{M^G(r, n+s)}) c(T_{Gr(\E_x, n)}) \over c(T_{Quot^n_C(\E)})
c(T_{Quot^n_{C'}(s)})} = {(c(Q) c(Q^*))^r (c(S^* \otimes Q)) \over
(c(Q))^r (c(Q))^r}=
$$
$$
= {(c(Q^*))^r (c(Q))^r \over (c(Q))^r c(Q \otimes Q^*)} =
{(c(Q^*))^r \over c(Q \otimes Q^*)} = c(Q^* \otimes S).
$$ 

  Hence the full Chern class $c(V)$ of the excess normal bundle $V$ is
  equal to the full Chern class $c(Q^* \otimes S)$ of the cotangent 
bundle to $Gr(\E_x, n)$ (which is isomorphic to $Q^* \otimes
  S$). Therefore
$$
c_{(r-n)n}(V) = (-1)^{rn-n^2} {r \choose n} = (-1)^{(r-1)n} {r \choose
  n}.
$$

 Recall that the top Chern class above  is equal to the intersection 
number of
$[Quot^n_C(\E)]$ and $[Quot^n_{C'}(s)]$. Now Theorem 4.7 follows from
the binomial formula. \quad \quad $\Box$

\section{Appendix: the punctual Quot scheme.}

In this appendix we give a proof of Theorem 2.2 saying that the
scheme $Quot(r, n)$ is irreducible of dimension $(rn-1)$.

Our strategy is   to  find a dense irreducible open subset  $W \subset 
\Q(r, n)$ of dimension $(rn-1)$. This subset will turn out to be a 
vector bundle over the punctual Hilber scheme $\Hi^n:= \Q(1, n)$
considered by Brian\c con and Iarrobino.

  We define $W$ as the set of all quotients $\Or 
\stackrel{\phi}\longrightarrow A$, $\phi=(\phi_1+ \phi_2 + \ldots + \phi_r)$
such that the first component $\phi_1: \Oo \to A$ is surjective (this
is clearly an open condition). Such a $\phi_1$ corresponds to a point
in $\Hi^n$. Once $\phi_1$ is chosen, the other components
$(\phi_2, \ldots, \phi_r)$ are given by an arbitrary element of
$Hom(\Oo^{\oplus(r-1)}, A) = \C^{(r-1)n}$. Therefore $W$ is a rank
$(r-1)n$ vector bundle over $\Hi^n$. By [Br] or [Ia] the subset 
 $W$ is irreducible of dimension $(rd-1)$. 

  Now we want to show that $W$ is dense in $\Q(r, n)$. In fact, for
any point $x \in \Q(r, n)$ we will find an irreducible rational curve
$C \subset \Q(r, n)$ connecting it with some point in $W$.

To that end,  we generalize Nakajima's  construction (cf. [Na, Chapter 2]) of  
the global Hilbert  scheme  $Hilb^n(\C^2)$ to the 
Quot scheme.   Once we  do that, the   existence  of the   
irreducible curve  will amount  to  an exercise in linear algebra 
(cf. Lemma 7.3).

\bigskip

 Fix a  complex  vector space $V$ of   dimension $n$, and  $N_n$ let 
be the space of pairs   of commuting nilpotent  operators  on $V$.  The space
$N_n$  is naturally a  closed    affine subvariety of   $\mathrm{End}(V)  
\oplus \mathrm{End}(V)$.    

  Consider a subspace $U_r$   of $N_n \times V^{\oplus r}$ formed 
by  all $(B_1,  B_2, v_1, \ldots,  v_r)$, such that there is 
no proper subspace of  $V$ which is invariant under $B_1, B_2$ and 
contains $v_1, \ldots, v_r$. Then $U_1 \times V^{\oplus (r-1)} \subset
U_2 \times V^{\oplus (r-2)} \subset \ldots \subset U_r$ is a chain
of open subsets in $N_n \times V^{\oplus r}$ (each of them is given by 
a condition saying that some system of vectors in $V$ has maximal rank).

\medskip

Ona has a natural $GL(V)$-action on $V_r$ and it is easy to prove that
$U_r$ is $GL(V)$-stable.

\begin{lemma} 
$GL(V)$ acts freely on $U_r$. 
\end{lemma}

{\it Proof.}   Suppose    $g \in  GL(V)$ stabilizes
$(B_1, B_2,   v_1, \ldots, v_r)   \in  U_r$. Then  $\mathrm{Ker}(1-g)$
contains $ v_1, \ldots, v_r$. Since it is also preserved by $B_1, B_2$
, we have $Ker(1-g)=V$ and therefore $g=1$. \hfill $\Box$

 The following lemma  gives  an explicit construction of  the
punctual Quot scheme:

\begin{lemma} There exists a morphism $\pi: U_r \to \Q(r, n)$ such that

(i) $\pi$ is surjective;

(ii) the fibers of $\pi$ are precisely the orbits of $GL(V)$ action on 
$U_r$;

(iii) $\pi^{-1}(W) = U_1 \times V^{\oplus (r-1)}$. 
\end{lemma}

{\it Proof.} Note that the punctual Quot scheme does not depend on 
the surface since the $n$-th power of the maximal ideal ${\frak m}_x$
acts by zero on any length $n$ sheaf supported at $x$. Hence
we can  assume that $S= \C^2  = Spec \;\C[x_1, x_2]$ and that all
the quotients are supported at $x=0  \in \C^2$.

To  construct $\pi$ suppose   that $(B_1, B_2, v_1, \ldots,
v_r)$  is a  point  in  $U_r$ and   consider  a $\C[x_1,  x_2]$-module
structure on $V$ in which $x_1$ acts by  $B_1$ and $x_2$ acts by $B_2$. 
We can view $V$ as a quotient of a free $\C[x_1,   x_2]$-module with 
generators  $v_1, \ldots,  v_r$.  Since $B_1$ and $B_2$ are  nilpotent 
$\sqrt{\mathrm{Ann}(V)} = (x_1, x_2)$. Therefore a  coherent sheaf $A$  
on  $\C^2$ associated  with $V$ is  a quotient of $\Or$ supported 
at $s$. Moreover,  $V \simeq H^0(S, A)$ as vector spaces.

A different point in the    same $GL(V)$-orbit defines an   isomorphic
quotient,  hence  the  fibers of $\pi:   U_r/  \mathrm{GL(V)} \to Quot
(r,n)$ are $GL(V)$-invariant.  Moreover, suppose that two points $u_1,
u_2$ of $U_r$ give rise to isomorphic quotients $A_1$, $A_2$. Then the
induced   isomorphism    between     $H^0(S,     A_1)$    and
$H^0(S, A_2)$ defines an  element of $GL(V)$ taking  $u_1$ to
$u_2$. Therefore, each fiber of $\pi$ is  precisely one $GL(V)$-orbit. 
This proves $(ii)$.

To  prove $(i)$, suppose we have  a quotient $ \Or  \to  A \to  0 $ of
length $d$ supported   at  zero. Multiplication  by   $x_1$ and  $x_2$
induces a pair   of commuting nilpotent  operators on $H^0(S,
A)$. Choose  a $\C$-linear isomorphism $H^0(S,  A) \simeq V$. 
The generators   of the free   $\C[x_1, x_2]$-module  $H^0(S,
\Or)$ project to some vectors $v_1,  \ldots, v_r$ in $V$.  Since $v_1,
\ldots, v_r$ generate $V$ as a $\C[x_1, x_2]$-module, $(x_1, x_2, v_1,
\ldots v_r)$ is a point of $U_r$. Thus $(i)$ is proved.

  Finally, $(iii)$ follows from definitions of $W$ and $U_1$.
\hfill $\Box$
 
\bigskip
 
  Now we want to show that any point in $U_r$ can be deformed
to a point in the preimage of $W$. The above construction will allow us to
construct this deformation using the following lemma

\begin{lemma}
Let $B_1$, $B_2$ be two commuting nilpotent operators on a vector space V.
There exists a third nilpotent operatov $B_2'$ and a vector $w \in V$
such that

(i) $B_2'$ commutes with $B_1$;

(ii) any linear combination $\alpha B_2 + \beta B_2'$ is nilpotent;

(iii) $(B_1, B_2', w) \in U_1$, i.e. $w$  is a cyclic vector for the
pair of operators $(B_1, B_2')$.
\end{lemma}

  This lemma will be proved later. Now we will show how it can be used
to give an

\bigskip 

\noindent\emph{End of Proof of Theorem 2.2:} 

Let $x$ be a point of $\Q(r, n)$ and $u_1 = (B_1, B_2, v_1, \ldots, v_r)$
be any point of $\pi^{-1}(x) \subset   U_r$. Choose a nilpotent
operator $B_2'$ and a vector $w \in V$ as in the lemma above. 
Connect the points $u_1$ and $u_2 = (B_1, B_2', w, v_2 , \ldots, v_r)$ with
a straight line $\Phi(t)$, $\;t \in \C$ such that  $\Phi(1) = u_1$ and 
$\Phi(0) = u_2$. This $\Phi(t)$ is given by equation:
$$
\Phi(t)= (B_1, tB_2' + (1-t)B_2, tw +(1-t)v_1, v_2, \ldots, v_r)
$$

  Note  that  for all $t \in C$,   $\; B_2(t) = tB_2' + (1-t)B_2 $   
is nilpotent and commutes with $B_1$. Therefore the image of $\Phi(t)$
is a subset of $N_n \times V^{\oplus r}$. Since $U_r$ is open in 
$N_n \times V^{\oplus r}$, there is a dense open subset 
$C \subset \C$ such that $\Phi(C) \subset U_r$.
Similarly, there exists a dense open subset $C_1 \subset 
C$ such that $\Phi(C_1) \subset U_1 \times V^{\oplus 
(r-1)}$.

  Hence the image $\pi(\Phi(C)) \subset \Q(r, n)$ is an 
irreducible rational curve connecting $x = \pi(u_1)$ with $\pi(u_2) \in W$.
Note that $\pi(\Phi(C_1)) \subset W$.
Therefore $x$ belongs to the closure of $W$. Since by Theorem 1.1
$W$ is irreducible of dimension $(rd-1)$, the scheme $\Q(r, n)$ is
also irreducible of dimension $(rd-1)$. Theorem 2.2 is proved. $\Box$

\bigskip 

\noindent\emph{Proof of Lemma 7.3:}

\medskip 

\emph{Step 1.}
 We will find a basis $e_{i,j}$ of $V$, where $1 \leq i \leq k$ and
$1 \leq j \leq \mu_i$ such that

\bigskip

(a) $B_1^{j-1}(e_{i, 1}) = e_{i, j}$ for $j \leq \mu_i$ and 
$B_1^{\mu_j}(e_{i, 1}) = 0$ (i.e. $B_1$ has Jordan canonical form
in the basis $e_{i, j}$);

(b) $B_2(e_{i, 1}) \in \Big(\bigoplus_{k \geq i +1} \C\cdot e_{k, 1} 
\Big) \oplus B_1 \cdot V$.

\bigskip 

 To that end, recall one way to construct a Jordan basis for $B_1$.
Let $d = \dim V$ and $V_i = Ker(B_1^{d-i})$. The subspaces $V_i$ form 
a decreasing filtration $V = V_0 \supset V_1 \supset V_2 \ldots $.
Moreover, $B_1 \cdot V_i \subset V_{i+1}$. Firstly, we choose a basis
$(w_1, \ldots, w_{a_1})$ of $W_1: = V_0/V_1$. Lift this basis to some 
vectors $e_{1, 1}, e_{2, 1}, \ldots, e_{a_1, 1}$ in $V_0$ and set
all $\mu_1, \ldots, \mu_{a_1}$ equal to $d$. Secondly, choose a basis 
$(w_{a_1 + 1}, \ldots, w_{a_2})$ of $W_2: = V_1/(B_1\cdot V_0 + V_2)$. 
Lift this basis to some  vectors $e_{a_1 + 1, 1}, e_{2, 1}, \ldots, 
e_{a_2, 1}$ in $V_1$ and set all $\mu_{a_1+1}, \ldots, \mu_{a_2}$ equal 
to $d-1$. Continue in this manner by choosing bases of the spaces
$W_{i+1} = V_i/ (B_1 \cdot V_{i-1} + V_{i+1})$ and lifting them to $V_i$.
This procedure gives us  vectors $e_{1, 1}, e_{2, 1}, \ldots,
e_{k, 1}$ and the formula (a) tells us how to define $e_{i, j}$ for
$j \geq 2$. It is easy to check that the system of vectors $\{e_{i, j}\}$
is in fact a basis of $V$.

\bigskip 

  If we want to have property (b)  we should be more careful with 
the choice of $w_i$.  Note that all the subspaces $V_i$ and $B_1 \cdot V_i$
are $B_2$-invariant. Therefore we have  an induced action of $B_2$ on
each $W_i$. We can choose our basis $(w_{a_{i-1} + 1}, \ldots, w_{a_i})$ 
of $W_i$ in such a way that $B_2(w_i) \in \bigoplus_{s=i+1}^{a_i} \C \cdot
 w_s$ for all $i \in \{a_{i-1} + 1, \ldots, a_i \}$.
This ensures that $(b)$ holds as well.

\bigskip

\emph{Step 2.} Define $B_2'$ by $B_2'(e_{i, j}) = e_{i+1, j}$ if 
$j \leq \mu_{i+1}$ and 0 otherwise. It is immediate that $B_2'$ is nilpotent
and that $[B_1, B_2']=0$. Let $w = e_{1, 1}$. Then $e_{i, j} = B_1^{j-1}
(B_2^{i-1}(w))$ hence $(B_1, B_2', w) \in U_1$.

\bigskip

\emph{Step 3.}
Note that both $B_2$ and $B_2'$ are lower-triangular with zeros on the 
diagonal in the basis of $V$ given by
$$
e_{1, 1}, e_{2, 1}, \ldots, e_{k, 1}, e_{1, 2}, e_{2, 2}, \ldots
 e_{k, 2}, \ldots
$$ 
Hence any linear combination of $B_2$ and $B_2'$ is also lower-triangular
and has zeros on the diagonal. Therefore $\alpha B_2 + \beta B_2'$ is
nilpotent for any complex $\alpha$ and $\beta$. The Lemma is proved. 
\hfill $\Box$

\end{document}